\documentclass[11pt]{amsart}
\usepackage{amsmath,amsthm, amscd, amssymb, amsfonts}
\usepackage[all]{xy}
\usepackage{mathrsfs}
\usepackage{enumitem}
\usepackage[T2A,T1]{fontenc}
\usepackage{hyperref}
\newcommand{\pinom}{\genfrac{[}{]}{0pt}{}}

\allowdisplaybreaks

\usepackage[ansinew]{inputenc}
\usepackage{graphicx,fancyhdr}

\usepackage[dvips, dvipsnames, usenames]{color}

\newcommand{\co}{\operatorname{co}}
\newcommand{\comment}[1]{}

\newcommand{\actizq}{\triangleright}
\newcommand{\actder}{\triangleleft}

\numberwithin{equation}{section}
\newtheorem{theorem}{Theorem}[section]
\newtheorem{lemma}[theorem]{Lemma}

\newtheorem{prop}[theorem]{Proposition}

\theoremstyle{definition}
\newtheorem{definition}[theorem]{Definition}
\newtheorem{example}[theorem]{Example}

\theoremstyle{remark}
\newtheorem{remark}[theorem]{Remark}

\newcommand{\pf}{\begin{proof}}
  \newcommand{\epf}{\end{proof}}

\newcommand{\bd}{ \mathbf{d}}
\newcommand{\kk}{ \mathbf{k}}
\newcommand{\h}{ \mathbf{h}}
\newcommand{\ku}{ \Bbbk}

\newcommand{\I}{\mathbb I}

\newcommand{\K}{\mathbb K}
\newcommand{\N}{\mathbb N}
\newcommand{\bP}{\mathbf{P}}

\newcommand{\Q}{\mathbb Q}

\newcommand{\V}{\mathbb V}

\newcommand{\Z}{\mathbb Z}

\renewcommand{\_}[1]{_{\left( #1 \right)}}
\renewcommand{\^}[1]{^{\left( #1 \right)}}

\newcommand{\hwt}{\widetilde{h}}
\newcommand{\hwtb}{\widetilde{\h}}

\newcommand{\cA}{\mathcal{A}}
\newcommand{\cB}{\mathcal{B}}

\newcommand{\Ss}{{\mathcal S}}

\newcommand{\hgo}{\mathfrak h}

\newcommand{\n}{\mathfrak n}

\newcommand{\id}{\operatorname{id}}

\newcommand{\ydh}{{}^{H}_{H}\mathcal{YD}}
\newcommand{\hyd}{\mathcal{YD}{}^{H}_{H}}

\newcommand{\ydV}{{}^{V_{\cB}^0}_{V_{\cB}^0}\mathcal{YD}}
\newcommand{\Vyd}{\mathcal{YD}{}^{V_{\cB}^0}_{V_{\cB}^0}}

\newcommand{\ot}{\otimes}

\newcounter{tabla}\stepcounter{tabla}

\title[Lusztig quantum divided power algebra]
{On the Hopf algebra structure of the Lusztig quantum divided power algebras}

\author[Nicol\'as Andruskiewitsch, Iv\'an Angiono and  Cristian Vay]
{Nicol\'as Andruskiewitsch, Iv\'an Angiono and  Cristian Vay}

\thanks{The work of N. A., I. A. and C. V. was partially supported by CONICET and Secyt (UNC). The work of C.V. was partially supported by Foncyt PICT 2016-3927}

\address{\noindent Facultad de Matem\'atica, Astronom\'{\i}a y F\'{\i}sica,
	Universidad Nacional de C\'ordoba. CIEM -- CONICET. 
	Medina Allende s/n (5000) Ciudad Universitaria, C\'ordoba,
	Argentina}
\email{(andrus|angiono|vay)@famaf.unc.edu.ar}

\subjclass[2010]{16T05}
\begin{document}

\begin{abstract}
We study the Hopf algebra structure of Lusztig's quantum groups.
First we show that the zero part is the tensor product of the group algebra of a finite abelian group
with the enveloping algebra of an abelian Lie algebra.
Second we build them from the plus, minus and zero parts by means of suitable actions and coactions within 
the formalism presented by Sommerhauser to describe triangular decompositions. 

\end{abstract}

\maketitle

\section{Introduction}\label{section:introduction}
  
There are two versions of quantum groups at roots of 1: the one introduced and studied by De Concini, Kac and Procesi \cite{DK-root of 1,DP-quantum}
and the  quantum divided power algebra of Lusztig \cite{L-AiM,L-contemporary,L-fdHa-JAMS,L-qgps-at-roots}. The small quantum groups 
(aka Frobenius-Lusztig kernels) appear as 
quotients of the first  and Hopf subalgebras of the second; in both cases they fit into suitable exact sequences of Hopf algebras. 

The key actor in all these constructions is what we now call a Nichols algebra of diagonal type.
Indeed all the Hopf algebras involved have triangular decompositions compatible with the mentioned exact sequences; 
the positive part of the small quantum group is a Nichols algebra.
The celebrated classification of the finite-dimensional Nichols algebras of diagonal type was achieved in 
\cite{H-classif RS}.  
The positive parts of the small quantum groups correspond to braidings of Cartan type, 
but there are also braidings of super and modular types in the list, see the survey \cite{AA17}.

The question of  defining the versions of the quantum groups of De Concini, Kac and Procesi on one hand, 
and of Lusztig on the other, for every Nichols algebra in the classification arises unsurprisingly.
The first was solved in \cite{A-preNichols} introducing Hopf algebras also with triangular decompositions and
whose positive parts are now the distinguished pre-Nichols algebras of diagonal type.
These were introduced earlier in \cite{Ang-crelle}, instrumental to the description of the defining relations of the Nichols algebras. A distinguished pre-Nichols algebra projects onto the corresponding Nichols algebra and the kernel
is a normal Hopf subalgebra that is even central under a mild technical hypothesis, see \cite{A-preNichols,AAR2}.
The geometry behind these new Hopf algebras is studied in \cite{AAY} for Nichols algebras coming in families.

Towards the second goal, the graded duals of those distinguished pre-Nichols algebras were studied in \cite{AAR1} under the name of Lusztig algebras; 
when the braiding is of Cartan type one recovers in this way the positive (and the negative) parts of  
Lusztig's quantum groups. 
A  Lusztig algebra contains the corresponding Nichols algebra as a normal Hopf subalgebra
and the cokernel is an enveloping algebra $U(\n)$ under the same mild technical hypothesis 
mentioned above, see \cite{AAR2}. 
In \cite{AAR2,AAR3} it was shown that $\n$ is either 0 or the positive part of a semisimple Lie algebra
that was determined explicitly in each case.

In order to construct the analogues of Lusztig's quantum groups at roots of one  for Nichols algebras of diagonal type, we still need to define the 0-part and the interactions with the positive and negative parts. 
This leads us to understand the Hopf algebra structure of a Lusztig's quantum group
 which is the objective of this Note.

Let $V$ be the $\Z[v, v^{-1}]$-Hopf algebra as in \cite[2.3]{L-fdHa-JAMS}; the quantum group is
defined by specialization of $V$. Our first goal is to describe the specialization of the 0-part $V^0$,
a commutative and cocommutative Hopf subalgebra of $V$.
We show in Theorem \ref{th:V0-isom}  that
it splits as the tensor product of the group algebra of a finite group
and the enveloping algebra of the Cartan subalgebra of the corresponding Lie algebra.
For this we use some skew-primitive elements $h_{i,n} \in V^0$, cf. Definition \ref{def:hn-skew-primitive},
defined from the elements $\pinom{K_i; 0}{t}$ and $K_i^{n}$ of the original presentation of \cite{L-fdHa-JAMS}. 
The elements $h_{i,n}$, or rather multiples of them, were already introduced in \cite{Lentner} towards 
defining unrolled versions of quantum groups; see Remark \ref{rem:lentner}.
We point out that the definition in \emph{loc. cit.} 
is by a limit procedure, while ours is explicit in terms of polynomials $p_{n,s}\in\Z[v,v^{-1}]$ that we define recursively in Lemma \ref{lem:primitive-linear-comb}. 
Theorem \ref{th:V0-isom}  appears in \cite{Lentner2,Lentner}.

In \cite{S} it is explained that Hopf algebras $U$ with a triangular decomposition $U \simeq A \ot H \ot B$,
where $H$ is a Hopf subalgebra, $A$ is a Hopf algebra in the category of left Yetter-Drinfeld modules and $B$ is the same but right, plus natural compatibilities, can be described by some specific structure that we call a TD-datum. 
Our second goal is to spell out the TD-datum corresponding to the quantum group, see Theorem \ref{th:structure-simply-laced}. 

The paper is organized as follows. In Section \ref{sec:preliminaries}
we set up some notation and recall the formalism of \cite{S}.
Section \ref{sec:V0} contains the analysis of the Hopf algebra $V^0$ from \cite{L-fdHa-JAMS}.
In Section \ref{sec:V-slaced} we recall the definition of Lusztig's version of quantum groups at roots of 1,
show that it fits into the setting of \cite{S} and prove Theorem \ref{th:structure-simply-laced}. For simplicity of the exposition we assume that the underlying  Dynkin diagram is simply-laced; in
the last Subsection we discuss how one would extend the material to the general case.

The Lusztig's quantum groups enter into a cleft short  exact sequence of Hopf algebras \cite[3.4.1,3.4.4]{A-canad} 
and contain an unrolled version of the finite quantum groups \cite{Lentner} 
but as is apparent from the description here, they are not unrolled Hopf algebras. 

We are not aware of other papers containing information on the matter of our interest.
Other versions of triangular decompositions similar to \cite{S} appear in \cite{Majid,Laugwitz}.

\section{Preliminaries}\label{sec:preliminaries}
  
\subsection{Conventions}\label{subsection:conventions}
We adhere to the notation in \cite{L-fdHa-JAMS,L-qgps-at-roots}
as much as possible.
If $t \in\N_0$, $\theta \in \N$ and $t < \theta$, then $\I_{t, \theta} := \{t, t+1,\dots,\theta\}$,
$\I_{\theta} := \I_{1, \theta}$. 

\smallbreak
Let $\cA = \Z[v, v^{-1}]$, the ring of Laurent polynomials in the indeterminate $v$,  
$\cA' = \Q (v)$, its field of fractions; later we also need $\cA'' := \Z[v, v^{-1}, (1-v^2)^{-1}]$. 
The $v$-numbers are the polynomials
\begin{align*}
[s]_{v} &=\frac{v^{s} - v^{-s}}{v - v^{-1}}, & [N]_{v}^!&=\prod_{s=1}^{N} [s]_{v}, &
\pinom{N}{i}_v & =\frac{[N]_v^!}{[N-i]_v^![i]_v^!} \in \cA,
\end{align*}
$s,N\in \N$, $i \in  \I_{0,N}$. We denote $\displaystyle\pinom{N}{i}_v = 0$ 
when $i > N$, $i,N \in \N$.


If $B$ is a commutative ring and $\xi \in B$ is a unit, then $B$ is an $\cA$-algebra 
via $v \mapsto \xi$; the elements $[s]_{v}$, $[N]_{v}^!$, $\pinom{N}{i}_v$
of $\cA$ specialize to $[s]_{\xi}$, $[N]_{\xi}^!$, $\pinom{N}{i}_{\xi}$ of $B$.
As in \cite[35.1.3]{L-libro}, we fix $\ell\in \N$ and set
\begin{align*}
\ell' = \begin{cases} \ell &\text{if $\ell$ is odd,}\\
2\ell &\text{if $\ell$ is even.} \end{cases}
\end{align*}
This convention is slighty different from the one in \cite[5.1, pp. 287 ff]{L-fdHa-JAMS}.

Let $\phi_{\ell'} \in \Z[v]$ be the $\ell'$-th cyclotomic 
polynomial; let  $\cB$ be the field of fractions of $\cA / \langle \phi_{\ell'} \rangle$
and let $\xi$ be the image of $v$ in $\cB$. We have in $\cB$
\begin{align}
\phi_{\ell} (\xi^2) &=0, & \xi^{\ell} &= (-1)^{\ell' + 1}, \qquad \xi^{2\ell} = \xi^{\ell^2} = 1,
\\ \label{eq:se-anula}
{\pinom{N + M}{M}}_{\xi} &= 0,& N, M &\in \I_{0,\ell - 1}, \qquad N+ M \geq \ell. 
\end{align}
We also have that
\begin{align}\label{eq:binom-evaluation-1}
\frac{[m\ell]_{\xi}}{[n\ell]_{\xi}} & = \xi^{(m-n)\ell}\frac{m}{n}, & 
\frac{[m\ell+j]_{\xi}}{[n\ell+j]_{\xi}} &= \xi^{(m-n)\ell}, & m,n\in\N_0, & \, j\in \I_{\ell - 1}.
\end{align}
Hence for all $m \ge n\in\N_0$ and $j\in \I_{\ell - 1}$, we have 
\begin{align}\label{eq:binom-evaluation-2}
{\pinom{m \ell}{n \ell}}_{\xi} &= \binom{m}{n}, & {\pinom{m \ell+j}{j}}_{\xi} &= \xi^{mj\ell}, & {\pinom{m \ell+j-1}{j}}_{\xi} &= 0.
\end{align}

Let $\ku$ be a field; all algebras, coalgebras, etc. below are over $\ku$ unless explicitly stated otherwise.
If $A$ is an associative unital $\ku$-algebra, then we identify $\ku$ with a subalgebra of $A$.

\subsection{Hopf algebras with triangular decomposition}\label{subsec:TD}

Let $H$ be a Hopf algebra with multiplication $m$, comultiplication $\Delta$ (with Sweedler notation $\Delta(h) = h\_{1} \otimes h\_{2}$), counit $\varepsilon$ 
and bijective antipode $\Ss$; we add a subscript $H$ when precision is desired. 
We denote by $\ydh$, respectively $\hyd$, the category of left-left, respectively right-right Yetter-Drinfeld modules
over $H$. 
If $M \in \ydh$, then the action of $H$ on $M$ is denoted by $\actizq$ 
while the coaction is $m \mapsto m\_{-1} \otimes m\_{0}$, 
whereas if $N \in \hyd$, the action is denoted by $\actder$ 
and the coaction is $n \mapsto n\_{0} \otimes n\_{1}$.
For Hopf algebras either in $\ydh$ or $\hyd$, we use notations as above but with the variation 
$\Delta(r) = r\^{1} \otimes r\^{2}$.
Given Hopf algebras $R$ in $\ydh$ and $S$ in $\hyd$, their bosonizations are denoted $R \# H$, $H \# S$.
If $A \overset{\pi}{\underset{\iota}{\rightleftarrows}} H$ are morphisms of Hopf algebras with $\pi\iota = \id_H$, then
$R \# H \simeq A \simeq H \# S$ where $R$ and $S$ are the subalgebras of right and left coinvariants of $\pi$; see \cite[\S 11.6, \S 11.7]{Rad-libro}.

\medbreak
A \emph{TD-datum} over $H$ \cite[Definition 3.2]{S} is a collection $(A, B, \rightharpoonup, \leftharpoonup, \sharp)$ where 

\begin{enumerate}[leftmargin=*,label=\rm{(\roman*)}]
	\item $A$ is a Hopf algebra in $\ydh$;
	\item $B$ is a Hopf algebra in $\hyd$;
	\item $\rightharpoonup: B \otimes A \to A$ is a left action so that $A$ is a left $B$-module,
and	the following identities hold for all $a \in A$, $b\in B$ and $h\in H$:
\begin{align}
&	\begin{aligned}
b \rightharpoonup (h \actizq a)
&= h\_{1} \actizq\left(\left(b \actder h\_{2}\right) \rightharpoonup a\right), \qquad \qquad b \rightharpoonup 1 =  \varepsilon_B(b),
\\
\Delta_{A}(b \rightharpoonup a) &= \left({b\^{1}}\_{0} \rightharpoonup a\^{1}\right) \otimes\left({b\^{1}}\_{1} \actizq\left(b\^{2} \rightharpoonup a\^{2}\right)\right);
\end{aligned}
\end{align}
	
	\item $\leftharpoonup: B \otimes A \to B$ is a right action so that $B$ is a right $A$-module,
	and	the following identities hold for all $a \in A$, $b\in B$ and $h\in H$:
\begin{align}
\begin{aligned}
(b \actder h) \leftharpoonup a
&= \left(b \leftharpoonup \left(h\_{1} \actizq a\right)\right) \actder h\_{2},
\qquad \qquad 1 \leftharpoonup a =  \varepsilon_A(a),
\\
\Delta_{B}(b \leftharpoonup a) &= \left(\left(b\^{1} \leftharpoonup a\^{1}\right) \actder {a\^{2}}\_{-1}\right) \otimes\left(b\^{2} \leftharpoonup {a\^{2}}\_{0}\right); 
\end{aligned}
\end{align}
both actions also satisfy for all $a \in A$ and $b\in B$:
\begin{align}
&\begin{aligned}
&\left(b\^{1} \rightharpoonup a\^{1}\right) \otimes\left(b\^{2} \leftharpoonup a\^{2}\right)
\\&= \left({b\^{1}}\_{1} \actizq\left(b\^{2} \rightharpoonup {a\^{2}}\_{0}\right)\right) \otimes\left(\left({b\^{1}}\_{0} \leftharpoonup a\^{1}\right) \actder {a\^{2}}\_{-1}\right)
\end{aligned}
\end{align}

	\item $\sharp: B \otimes A \to H$ is a linear map, $b \otimes a \mapsto b \sharp a$,
satisfying the following identities for all $a, c \in A$, $b,d\in B$, $h \in H$.
\end{enumerate}

\medbreak
\noindent \emph{Compatibility of $\sharp$ with the structure of $H$:}
\begin{align}
\begin{aligned}
\left(b \sharp\left(h\_{1} \actizq a\right)\right) h\_{2}
&= h\_{1}\left(\left(b \actder h\_{2}\right) \sharp a\right),
\\
\Delta_{H}(b \sharp a) &=\left({b\^{1}}\_{0} \sharp a\^{1}\right) {a\^{2}}\_{-1} \otimes {b\^{1}}\_{1}\left(b\^{2} \sharp {a\^{2}}\_{0}\right),\\
\varepsilon_H(b \sharp a) &= \varepsilon_B(b) \varepsilon_A(a),
\end{aligned}
\end{align}

\medbreak
\noindent \emph{Compatibility of $\sharp$ with the products of $A$ and $B$:}
\begin{align}
\begin{aligned}\label{eq:Compatibility sharp}
b \sharp\left(a c\right) &= \left(b\^{1} \sharp a\^{1}\right) {a\^{2}}\_{-1}\left(\left(b\^{2} \leftharpoonup {a\^{2}}\_{0}\right) \sharp c\right),
\\
\left(b d\right) \sharp a &= \left(b \sharp\left({d\^{1}}\_{0} \rightharpoonup a\^{1}\right)\right) {d\^{1}}\_{1}\left(d\^{2} \sharp a\^{2}\right),
\\
b \sharp 1 &=  \varepsilon_B(b), \qquad 1 \sharp a =  \varepsilon_A(a),
\end{aligned}
\end{align}

\medbreak
\noindent \emph{Compatibility of the actions with the multiplications via $\sharp$:}
\begin{align}
&\begin{aligned}\label{eq:Compatibility harpoonup}
b \rightharpoonup\left(a c\right) &= \left({b\^{1}}\_{0} \rightharpoonup a\^{1}\right) \times
\\ 
&\times
\big({b\^{1}}\_{1}(b\^{2} \sharp a\^{2}) {a\^{3}}\_{-1} \actizq \big[\big(b\_{3} \leftharpoonup {a\^{3}}\_{0}\big) \rightharpoonup c\big]\big),
\\
\left(b d\right) \leftharpoonup a &= \left(\left[b \leftharpoonup\left({d\^{1}}\_{0} \rightharpoonup a\^{1}\right)\right] \actder {d\^{1}}\_{1}\big(d\^{2} \sharp a\^{2}\big) {a\^{3}}\_{-1}\right) \times
\\
&\times \left(d\^{3} \leftharpoonup {a\^{3}}\_{0}\right);
\end{aligned}
\end{align}

\medbreak
\noindent \emph{Compatibility of the coactions with the comultiplications via $\sharp$:}
\begin{align}
&\begin{aligned}
&\left({b\^{1}}\_{0} \rightharpoonup a\^{1}\right)\_{-1} {b\^{1}}\_{1}\left(b\^{2} \sharp a\^{2}\right) \otimes\left({b\^{1}}\_{0} \rightharpoonup a\^{1}\right)\_{0}
\\&= \left({b\^{1}}\_{0} \sharp a\^{1}\right) {a\^{2}}\_{-1} \otimes\left({b\^{1}}\_{1} \actizq\left(b\^{2} \rightharpoonup {a\^{2}}\_{0}\right)\right);
\\
&\left(b\^{2} \leftharpoonup {a\^{2}}\_{0}\right)\_{0} \otimes\left(b\^{1} \sharp a\^{1}\right) {a\^{2}}\_{-1}\left(b\^{2} \leftharpoonup {a\^{2}}\_{0}\right)\_{1}
\\&= \left(\left({b\^{1}}\_{0} \leftharpoonup a\^{1}\right) \actder {a\^{2}}\_{-1}\right) \otimes {b\^{1}}\_{1}\left(b\^{2} \sharp {a\^{2}}\_{0}\right).
\end{aligned}
\end{align}

\begin{prop}\label{prop:sommerh}
\begin{enumerate}[leftmargin=*,label=\rm{(\alph*)}]
\item\label{item:prop-sommerh-1} \cite[3.3, 3.4]{S} 
Let $(A, B, \rightharpoonup, \leftharpoonup, \sharp)$ be a
TD-datum over $H$. Then $U :=  A \otimes H \otimes B$
is a Hopf algebra with multiplication, comultiplication and antipode:
\begin{multline*}
(a \otimes h \otimes b) (c \otimes k \otimes d) 
= a\left(h\_{1} \actizq \left({b\^{1}}\_{0} \rightharpoonup c\^{1}\right)\right) 
\\ 
\otimes h\_{2} {b\^{1}}\_{1}\left(b\^{2} \# c\^{2}\right) {c\^{3}}\_{1} k\_{1} \otimes\left(\left(b\^{3} \leftharpoonup {c\^{3}}\_{2}\right) \actder k\_{2} \right) d,
\\
\Delta(a \otimes h \otimes b)  = \left(a\^{1} \otimes {a\^{2}}\_{-1} h\_{1} \otimes {b\^{1}}\_{0}\right) \otimes\left({a\^{2}}\_{0} \otimes h\_{2} {b\^{1}}\_{1} \otimes b\^{2}\right),
\\
\Ss(a \otimes h \otimes b)=(1 \otimes 1 \otimes \Ss_{B}(b\_{0}))(1 \otimes \Ss_{H}(a\_{-1} h b\_{1}) \otimes 1)(\Ss_{A}(a\_{0}) \otimes 1 \otimes 1).
\end{multline*}

\item\label{item:prop-sommerh-2} \cite[3.5]{S} 
	Let $\mathcal U$ be a Hopf algebra. 
Let $A$ and $B$ be Hopf algebras in $\ydh$  and $\hyd$ respectively, provided with injective algebra maps
\begin{align*}
\iota_A: A &\hookrightarrow\mathcal U,& \iota_H: H &\hookrightarrow\mathcal U,& \iota_B: B &\hookrightarrow\mathcal U.
\end{align*}
 Assume that 
\end{enumerate}
\begin{enumerate}[leftmargin=*,label=\rm{(\roman*)}]
	\item The map $\xymatrix{A \ot H \otimes B \ar@{->}[rrr]^-{m_U(\iota_A\otimes \iota_H \otimes \iota_B)} & & & \mathcal U}$  
is a linear isomorphism.
	\item The induced maps $A \# H \to\mathcal U$, $H \# B \to\mathcal U$ are Hopf algebra maps.
\end{enumerate}
Then there exists a TD-datum $(A, B, \rightharpoonup, \leftharpoonup, \sharp)$ 
over $H$  such that $\mathcal U \simeq U$. \qed
\end{prop}

Clearly these constructions are mutually inverse.
In the setting of the Proposition, we say that $U \simeq A \ot H \otimes B$ is a \emph{triangular decomposition}.

As observed in \cite{S}, the verification of the conditions in the definition of TD-datum
is easier when $H$ is commutative and cocommutative.

\section{The algebra $V^0$}\label{sec:V0}
\subsection{Basic definitions}
We fix $\theta\in\N$. For simplicity, set $\I:=\I_{\theta}$. 
Following \cite[2.3, pp. 268 ff]{L-fdHa-JAMS} we consider the $\cA$-algebra $V^0$ presented by generators
\begin{align}\label{eq:V0-gens}
&K_i, & &K_i^{-1}, & &\pinom{K_i; c}{t}, & i\in \I, \, c &\in \Z, \, t\in \N_0
\end{align}
and relations for all $i\in \I$, tagged as in \emph{loc. cit.},
\begin{align}
\tag{g5} (v -v^{-1}) &\pinom{K_i; 0}{1} = K_i - K_i^{-1}, \label{eq:g5}
\\ \tag{g6}  &\text{the generators \eqref{eq:V0-gens} commute with each other,} \label{eq:g6}
\\ \tag{g7} K_i K_i^{-1} & = 1, \quad \pinom{K_i; 0}{0} = 1, \label{eq:g7}
\\ \tag{g8} \label{eq:g8}
{\pinom{t+t'}{t}}_{v} &\pinom{K_i; 0}{t + t'} = 
\sum_{0 \le j \le t'} (-1)^j v^{t(t'-j)} {\pinom{t+j-1}{j}}_{v} K_i^j \pinom{K_i; 0}{t}\pinom{K_i; 0}{t' - j},
\\ &\notag \hspace{180pt} t \geq 1, \, t'\geq 0,
\\ \tag{g9} \pinom{K_i; -c}{t}& = 
\sum_{0 \le j \le t} (-1)^j v^{c(t-j)} {\pinom{c+j-1}{j}}_{v} K_i^j \pinom{K_i; 0}{t - j},
\quad t \geq 0, \, c\geq 1, \label{eq:g9}
\\ \tag{g10} \pinom{K_i; c}{t}& = 
\sum_{0 \le j \le t}  v^{c(t-j)} {\pinom{c}{j}}_{v} K_i^{-j} \pinom{K_i; 0}{t - j},
\quad t \geq 0, \, c\geq 0. \label{eq:g10}
\end{align}

Observe that \eqref{eq:g9} and \eqref{eq:g10} actually define the elements $\displaystyle \pinom{K_i; c}{t}$, $c\in \Z - 0$, 
in terms of $K_i^{\pm 1}$ and 
\begin{align}\label{eq:kit}
k_{i,t} &:= \pinom{K_i; 0}{t}, & t\in \N, \, i&\in\I.
\end{align}
See also \S \ref{subsec:V-mlaced} for an equivalent formulation.
Set
\begin{align}\label{eq:ait}
a_{i,t}&= \frac{v^{-t} K_i - v^tK_i^{-1}}{v -v^{-1}}&
i\in \I, \ t &\in \Z.
\end{align}
Thus $\Ss(a_{i,t}) = - a_{i,-t}$.
Taking $t' = 1$ in \eqref{eq:g8} we have 
\begin{align}\label{eq:kt0}
 [t+1]_{v} k_{i,t+1} &= k_{i,t} (v^t k_{i,1} - [t]_{v}K_i) = k_{i,t} a_{i,t},
\\ \intertext{hence}
\label{eq:kt}
[t]_{v}^! k_{i,t}&= \prod_{0 \le s < t} a_{i,s}.
\end{align}
Multiplying \eqref{eq:kt0} by $v - v^{-1}$, we get
\begin{align}\label{eq:prop214-lu}
K_i^{2} k_{i,t} &= v^t (v^{t + 1}  - v^{-t - 1}) K_i k_{i,t+1} + v^{2t} k_{i,t}.
\end{align}

\begin{prop}\label{prop:V0-prop-lusztig}  
\begin{enumerate}[leftmargin=*,label=\rm{(\alph*)}]
\item\label{item:V0-prop-lusztig-1} \cite[Lemma 2.21]{L-fdHa-JAMS} The $\cA$-module $V^0$ is free with basis 
\begin{align}
\label{eq:V0-base-lusztig} &K_1^{\delta_1} \cdots K_\theta^{\delta_\theta} k_{1,t_1} \cdots k_{\theta,t_\theta}, & \delta_i &\in \{0,1\}, \, t_i\in \N_0, \, i\in\I.
\end{align}

\item\label{item:V0-prop-lusztig-2} \cite[2.22]{L-fdHa-JAMS} $V^0 \otimes_\cA \cA' \simeq \cA' [\Z^{\I}]$ as $\cA'$-algebras. \qed
\end{enumerate}
\end{prop}

Thus $V^0$ is an $\cA$-form of the group algebra $\cA' [\Z^{\I}]$; actually it is a form of the Hopf algebra structure
as we see next.

\begin{lemma}\label{lem:comult-V0}
The $\cA$-algebra $V^0$ is a Hopf algebra with comultiplication determined by
\begin{align}
\label{eq:comult-V-3} \Delta (K_i^{\pm 1})  &=  K_i^{\pm 1}  \otimes K_i^{\pm 1}&
i&\in \I.
\end{align}
\end{lemma}
\pf Since $V^0$ is a subalgebra of the  Hopf algebra $\cA' [\Z]$, we need to see that $\Delta(V^0) \subset V^0 \otimes_\cA V^0$. 
By \eqref{eq:g9} and \eqref{eq:g10}, it is enough to show that
\begin{align}
\label{eq:comult-V-4} \Delta (k_{i,t})  &= \sum_{0 \le s \le t} k_{i,t-s} K_i^{-s}  \otimes k_{i,s} K_i^{t-s} 
\end{align}
for all $i\in \I$ and $t \in \N$.
We proceed by induction on $t$. If $t=1$, then
\begin{align*}
\Delta(k_{i,1}) = \frac{1}{v - v ^{-1}} \left(K_i\otimes K_i - K_i^{-1} \otimes K_i^{-1} \right) = k_{i,1} \otimes K_i + K_i^{-1} \otimes k_{i,1}.
\end{align*}
If \eqref{eq:comult-V-4} is valid for $t$, then
\begin{align*}
\Delta & (k_{i,t+1}) = \frac{1}{[t+1]_{v}} \Delta(k_{i,t})\Delta(a_{i,t}) 
= \frac{1}{[t+1]_{v}} \Big(\sum_{0 \le s \le t} k_{i,t-s} K_i^{-s}  \otimes k_{i,s} K_i^{t-s} \Big) \times
\\ & \qquad \qquad \left( \frac{v^{-t} K_i \ot K_i - v^t K_i^{-1} \ot K_i^{-1}}{v -v^{-1}} \right)
\\ &= \sum_{0 \le s \le t} \frac{v^{-t} k_{i,t-s} K^{1-s}  \otimes k_{i,s} K^{t+1-s} 
- v^t k_{i,t-s} K_i^{-1-s}  \otimes k_{i,s} K_i^{t-1-s}}
{[t+1]_{v}(v -v^{-1})}
\\ &= 
\sum_{0 \le s \le t} \frac{v^{-s}}{[t+1]_{v}} \left([t-s+1]_{v} k_{t+1-s} + 
\frac{v^{t-s} k_{i,t-s}K_i^{-1}}{v -v^{-1}} \right) K_i^{-s}  \otimes k_{i,s} K_i^{t+1-s}  
\\ & \quad + \frac{1}{[t+1]_{v}} \sum_{0 \le s \le t} v^{t-s} k_{i,t-s} K_i^{-1-s}  \otimes \left( [s+1]_{v} k_{i,s+1}  - \frac{v^{-s}k_{i,s} K_i}{v -v^{-1}} \right)  K_i^{t-s}
\\ &= 
k_{i,t+1}\otimes K_i^{t+1}+ \frac{1}{[t+1]_{v}}
\sum_{1 \le s \le t} v^{-s} [t-s+1]_{v} k_{i,t+1-s} K_i^{-s}  \otimes k_{i,s} K_i^{t+1-s}  
\\ & \qquad +  \sum_{0 \le s \le t-1} v^{t-s} \frac{[s+1]_{v}}{[t+1]_{v}} k_{i,t-s} K_i^{-1-s}  \otimes  k_{i,s+1} K_i^{t-s}
+ K_i^{-1-t}  \otimes k_{i,t+1}
\\ &= 
k_{i,t+1}\otimes K_i^{t+1} + K_i^{-1-t}  \otimes k_{i,t+1} 
\\ & \qquad + \sum_{1 \le j \le t} \frac{v^{-j} [t-j+1]_{v} + v^{t+1-j} [j]_{v}}{[t+1]_{v}}  k_{i,t+1-j} K_i^{-j}  \otimes k_{i,j} K_i^{t+1-j} 
\\ & = \sum_{0 \le s \le t+1} k_{i,t+1-s} K_i^{-s}  \otimes k_{i,s} K_i^{t+1-s},
\end{align*}
which completes the inductive step.
\epf 

\subsection{Some skew-primitive elements}

We introduce some notation:

\begin{itemize}[leftmargin=*]
\item For $n\in\N$, $\phi_n:\N\to\{0,1\}$ is the map given by $\phi_n(j)=0$ if $n-j$ is even and $\phi_n(j)=1$ if $n-j$ is odd.
\item $\Phi:V^0\to V^0$ is the algebra automorphism determined by
\begin{align}\label{eq:autom-1}
K_i & \mapsto -K_i, & K_i^{-1} & \mapsto -K_i^{-1}, & 
\pinom{K_i; c}{t} & \mapsto (-1)^t \pinom{K_i; c}{t}.
\end{align}
\end{itemize}
It is easy to see that \eqref{eq:autom-1} defines an algebra map.
Notice that 
\begin{align}\label{eq:autom-1-props}
\Phi(k_{i,n})&=(-1)^n k_{i,n},& \Phi(K_i^{\pm n})&=(-1)^n 
K_i^{\pm n}, & n &\in \N, \ i\in \I.
\end{align}

\begin{lemma}\label{lem:primitive-linear-comb} Let $n \in \N$.
We define $p_{n,s}\in\Z[v,v^{-1}]$, $s\in\I_n$, recursively on $s$ by 
$p_{n,1}= v^{-\phi_n(1)}$,
\begin{align*}
p_{n,s} &=\frac{v^{ns}-v^{-ns}}{v^{\phi_n(s)s}(v^n-v^{-n})} - \sum_{t \in \I_{s-1}} p_{n,t} {\pinom{s}{t}}_{v} v^{(\phi_n(t)-\phi_n(s))s}, &  &s>1.
\end{align*}
Then
\begin{align}\label{eq:Kn-K-n-equality}
K_i^{n}-K_i^{-n} &= (v^n-v^{-n}) \sum_{s \in \I_n} p_{n,s} \, k_{i,s} K_i^{\phi_n(s)}.
\end{align}

\end{lemma}

\noindent \emph{Proof.} Fix $i\in\I$.
By Proposition \ref{prop:V0-prop-lusztig} \ref{item:V0-prop-lusztig-1}, $K_i^{n}-K_i^{-n}$ is a linear combination of $k_{i,t}$, $k_{i,t} K_i$, $t\in\N_0$.
Indeed, it can be shown by induction on $n$ that $K_i^{\pm n}$ belongs to the $\cA$-submodule spanned by $k_{i,t}$, $k_{i,t} K_i$, $t\le n$. 
Using the involution $\Phi$, we see by \eqref{eq:autom-1-props}  
 that there are $a_{n,t} \in\cA$, $t\in\I_n$ such that
\begin{align}\label{eq:Kalan-linearcomb-kj}
K_i^{n}-K_i^{-n} &= \sum_{t \in \I_{0,n}}a_{n,t} k_{i,t} K_i^{\phi_n(t)}.
\end{align}

We extend scalars as in Proposition \ref{prop:V0-prop-lusztig} \ref{item:V0-prop-lusztig-2} and consider the algebra maps
\begin{align}\label{eq:defn-xi-character}
\Xi_{i,j}: &\cA' [\Z^{\I}] \to \cA', &  &K_i\mapsto v^j, & K_p &\mapsto 1, \quad p\neq i,
\end{align}
$j \in \N_{0}$. Notice that, with the convention  $\pinom{N}{n}_v = 0$ when $n > N$,
\begin{align}\label{eq:maps-preserve}
\Xi_{i,j}(k_{i,t}) &= 
\frac{1}{[t]_{v}^!} \prod_{0 \le s < t} \Xi_{i,j}(a_{i,s}) = 
\frac{1}{[t]_{v}^!} \prod_{0 \le s < t} \frac{v^{j-s} - v^{s-j}}{v -v^{-1}}
=  {\pinom{j}{t}}_{v}.
\end{align}
Applying $\Xi_{i,0}$ and $\Xi_{i,1}$ to \eqref{eq:Kalan-linearcomb-kj}, we see that
$0=a_{n,0}$, $v^n-v^{-n}=a_{n,1}v^{\phi_n(1)}$. Now we apply $\Xi_{i,s}$, $s>1$, to \eqref{eq:Kalan-linearcomb-kj}:
\begin{align*}
v^{ns} &- v^{-ns} = \sum_{t \in \I_{s}} a_{n,t} {\pinom{s}{t}}_{v} v^{\phi_n(t)s};
\\ \intertext{this implies the recursive formula holds since}
a_{n,s}  &= v^{-\phi_n(s)s}(v^{ns}-v^{-ns}) - \sum_{t \in \I_{s-1}} a_{n,t} {\pinom{s}{t}}_{v} v^{(\phi_n(t)-\phi_n(s))s}.
\qed
\end{align*}

\begin{definition}\label{def:hn-skew-primitive}
Let $n\in\N$. We set 
\begin{align}\label{eq:defn-hn}
h_{i,n} &:= \frac{K_i^{n}-K_i^{-n}}{n(v^n-v^{-n})} K_i^{n} = \frac{1}{n}\Big(\sum_{s \in \I_n} p_{n,s} \, k_{i,s} K_i^{\phi_n(s)} \Big) K_i^{n} \in V^0.
\end{align} 
Then $h_{i,n} = \dfrac{K_i^{2n}-1}{n(v^n-v^{-n})}$ is $(1,K_i^{2n})$-skew primitive.
\end{definition}

\begin{remark}\label{rem:lentner}
The elements $H'_\alpha$ defined in \cite[Theorem 3.1]{Lentner} are multiples of the above elements in the particular case $n=\ell$. Explicitly, $H'_{\alpha}=\frac{\ell(v^{\ell}-v^{-\ell})}{\phi_{\ell'}(v^{2})}\h_{i,\ell}$. Notice that $H'_{\alpha}$ are defined by taking a limit while \eqref{eq:defn-hn} is an explicit expression in terms of the polynomials $ p_{n,s}$ that are defined recursively. We discuss now these polynomials.
\end{remark}

\begin{lemma}\label{rem:pnn} Let $n\in\N$. Then 
\begin{align}\label{eq:pnn}
p_{n,n} =v^{-\binom{n}{2}} (-1)^{n-1}(v-v^{-1})^{n-1}[n-1]_{v}^!.
\end{align}
\end{lemma}

\pf We compute $k_{i,t}$ in $V^0 \otimes_\cA \cA' \simeq \cA' [\Z^\I]$:
\begin{align*}
k_{i,t} &= \frac{1}{[t]_{v}^!} \prod_{j=0}^{t-1} \frac{v^{-j}K_i-v^jK_i^{-1}}{v-v^{-1}} = \sum_{s=-t}^t f_{t,s} K_i^s,
& \text{for some }f_{t,s}&\in \cA'.  
\end{align*}
In particular, $f_{t,-t} =\dfrac{(-1)^{t}v^{\binom{t}{2}}} {[t]_{v}^!(v-v^{-1})^t}$.
Looking at the equality \eqref{eq:Kn-K-n-equality}, $K_i^{-n}$ appears only in one summand, $p_{n,n}k_{i,n}$, on the right hand side. Hence 
\begin{align*}
-1 &= (v^n-v^{-n})f_{n,-n} p_{n,n} = [n]_{v}(v-v^{-1}) \frac{(-1)^{n}v^{\binom{n}{2}}} 
{[n]_{v}^!(v-v^{-1})^n} p_{n,n},
\end{align*}
and the claimed equality  follows.
\epf

Given $\ell \in \N$ we consider the lower triangular matrix  
\begin{align*}
\bP_{i, \ell} &= \begin{pmatrix}
p_{11} & 0 & 0 &\dots &\dots &0 \\
p_{21}K_i & p_{22} & 0 &\dots &\dots &0 \\
p_{31} & p_{32}K_i & p_{32} &\dots &\dots &0 \\
\vdots & \vdots & \vdots & \vdots & \vdots & \vdots 
\\ p_{\ell 1} K_i^{\phi_\ell(1)} & p_{\ell 2}K_i^{\phi_\ell(2)} & p_{\ell 3} K_i^{\phi_\ell(3)}
&\dots & p_{\ell \ell-1}K_i & p_{\ell \ell} \end{pmatrix},&
\end{align*}
and the column vectors 
\begin{align*}
\kk_{i, \ell} &= \begin{pmatrix}
k_{i, 1} \\ \vdots \\ k_{i, \ell}
\end{pmatrix},& \hwtb_{i, \ell} &= \begin{pmatrix}
\hwt_{i, 1} \\ \vdots \\ \hwt_{i, \ell}
\end{pmatrix},& &\text{where } \hwt_{i,n} := nh_{i,n} K_i^{-n}.
\end{align*}
Then \eqref{eq:defn-hn} says that $\bP_{i, \ell}\kk_{i, \ell} = \hwtb_{i, \ell}$.
Recall $\cA'' = \Z[v, v^{-1}, (1-v)^{-1}]$ so that the matrix $\bP_{i, \ell}$ becomes invertible
in $V^0 \otimes_\cA \cA''$ by \eqref{eq:pnn}. Let
\begin{align*}
\bP_{i, \ell}^{-1} &= \begin{pmatrix}
q_{11} & 0 & 0 &\dots &\dots &0 \\
q_{21} & q_{22} & 0 &\dots &\dots &0 \\
q_{31} & q_{32} & q_{32} &\dots &\dots &0 \\
\vdots & \vdots & \vdots & \vdots & \vdots & \vdots 
\\ q_{\ell 1}  & q_{\ell 2} & q_{\ell 3} 
&\dots & q_{\ell \ell-1} & q_{\ell \ell} \end{pmatrix}.
\\
\intertext{Then}
k_{i,n} & =  \sum_{s \in \I_n} q_{n,s} \, \hwt_{i,s} =   \sum_{s \in \I_n} q_{n,s} s \, h_{i,s} K_{i}^{-s}, & n &\in \N, \ i\in \I.
\end{align*} 

\begin{example}
We compute $p_{n,s}$ for small values of $n$. For $n=2$, $p_{2,1}=v^{-1}$,
\begin{align*}
p_{2,2} &=\frac{v^{4}-v^{-4}}{v^2-v^{-2}} - p_{2,1} {\pinom{2}{1}}_{v} v^{2}
= v^2+v^{-2}-v(v+v^{-1})=v^{-2}-1.
\end{align*}
This agrees with \eqref{eq:pnn}. Thus,
\begin{align*}
h_{i,2}  &= \Big(\frac{v^{-2}-1}{2} \, k_{i,2} + \frac{v^{-1}}{2} \, k_{i,1} K_i\Big) K_i^2.
\end{align*}

For $n=3$, we have that $p_{3,1}=1$,
\begin{align*}
p_{3,2} &=\frac{v^{6}-v^{-6}}{v^{2}(v^3-v^{-3})} - p_{3,1} {\pinom{2}{1}}_{v} v^{-2}
= \frac{(v-v^{-1})^2[2]_{v}}{v^2},
\\
p_{3,3} &=\frac{v^{9}-v^{-9}}{v^3-v^{-3}} 
- \sum_{t \in \I_{2}} p_{3,t} {\pinom{3}{t}}_{v} v^{3\phi_3(t)}
= \frac{v^{9}-v^{-9}}{v^3-v^{-3}}
-p_{3,1} {\pinom{3}{1}}_{v}
-p_{3,2} {\pinom{3}{2}}_{v} v^{3}
\\ &
= 1-v^{-2}-v^{-4}+v^{-6} 
= \frac{(v-v^{-1})^2[2]_{v}}{v^3}.
\end{align*}
Again this agrees with \eqref{eq:pnn}. Thus,
\begin{align*}
h_{i,3}  &=  \big(\frac{(v-v^{-1})^2[2]_{v}}{3v^3} k_{i,3}  +  \frac{(v-v^{-1})^2[2]_{v}}{3v^2} k_{i,2} K_i + k_{i,1}\big) K_i^{3}.
\end{align*}
\end{example}

The element $H'$ computed in \cite[Example 3.2]{Lentner} (assuming  $\ell'=4$ in our notation)
is a multiple of $h_{i,2}$ above.

\begin{remark}  For instance,
from the preceding formulas we conclude:
\begin{align}\label{eq:ki1}
k_{i,1} &= h_{i,1}K_{i}^{-1},
\\\label{eq:ki2}
k_{i,2} &=  \frac{2}{v^{-2}-1}h_{i,2}K_{i}^{-2} - \frac{v^{-1}}{v^{-2}-1} h_{i,1},
\\\label{eq:ki3}
k_{i,3} &= \frac{3v^3}{(v-v^{-1})^2[2]_{v}} h_{i,3}K_i^{-3}
-\frac{2v}{v^{-2}-1}h_{i,2}K_{i}^{-1} + \frac{1}{v^{-2}-1} h_{i,1}K_i
\\ \notag &-\frac{3v^3}{(v-v^{-1})^2[2]_{v}} h_{i,1}K_i^{-1}.
\end{align}
\end{remark}

\subsection{Specializations of $V^0$}
Recall that $\ell' \in \N$ is defined in  Section \ref{subsection:conventions} and that
 $\cB$ is the field of fractions of 
$ \cA / \langle \phi_{\ell'} \rangle$.
We study now $V^0_{\cB} := V^0 \otimes_\cA \cB$. 
Thus the map  $\cA \to \cB$ factorizes through $\cA''$.

\begin{lemma}\cite[Lemma 4.4]{L-contemporary}, \cite[Lemma 2.21]{L-fdHa-JAMS}\label{lemma:V0-generators}
 The algebra $V^0_{\cB}$ is generated by $K_i$ and 
$\displaystyle k_{i,\ell} = \pinom{K_i;0}{\ell}$, $i\in\I$. Furthermore,
\begin{align}
K_i^{2\ell}=1.\label{eq:K2l=1}
\end{align}
\end{lemma}

\pf We first prove \eqref{eq:K2l=1}.
Taking $t=\ell-1$ and $t' = 1$ in \eqref{eq:g8} we have:
\begin{align*}
0 &\overset{\eqref{eq:se-anula}}{=} {\pinom{\ell}{\ell-1}}_{\xi} \pinom{K_i;0}{\ell} = \xi^{\ell-1}\pinom{K_i;0}{\ell-1}\pinom{K_i;0}{1} - {\pinom{\ell-1}{1}}_{\xi} K_i \pinom{K_i;0}{\ell-1}
 \\ &
= \pinom{K_i;0}{\ell-1} \frac{\xi K_i-\xi^{-1}K_i^{-1}}{\xi-\xi^{-1}} 
\overset{\eqref{eq:kt}}{=} \frac{\prod_{0 \le s < \ell}\xi^{-s}K_i-\xi^s K_i^{-1}}{[\ell - 1]_{\xi}^!(\xi-\xi^{-1})^{\ell}}
 \\ &=\xi^{-\binom{\ell}{2}} \frac{K_i^{\ell}-K_i^{-\ell}}{[\ell - 1]_{\xi}^! (\xi-\xi^{-1})^{\ell}}.
\end{align*}
Hence \eqref{eq:K2l=1} holds.

Let $t\in\I_{0,\ell - 1}$. Then $[t]_{\xi}^!\ne 0$, so 
$k_{i,t}= \frac{1}{[t]_{\xi}^!} \prod_{0 \le s < t} a_{i,s}$ belongs to the subalgebra generated by $K_i$.
Now we claim that
\begin{align}\label{eq:K0nl-formula}
k_{i,n\ell} &= \frac{1}{n!}\prod_{0 \le s < n} \left(k_{i,\ell}-s K_i^{\ell}\right), & \text{for all }&n\in\N.
\end{align}
We take $t=n\ell$ and $t' = \ell$ in \eqref{eq:g8}. By \eqref{eq:binom-evaluation-2} the left-hand side is
\begin{align*}
{\pinom{(n+1)\ell}{n\ell}}_{\xi} &k_{i, (n+1)\ell} 
= (n+1) k_{i, (n+1)\ell},
\end{align*}
while the right-hand side is, by \eqref{eq:binom-evaluation-1},
\begin{align*}
\sum_{0 \le j \le \ell} & (-1)^j \xi^{n\ell(\ell-j)} {\pinom{n\ell+j-1}{j}}_{\xi} K_i^j k_{i,n\ell}k_{i,\ell-j} 
\\ & = k_{i,n\ell}k_{i,\ell} 
+(-1)^{\ell}\xi^{n\ell(\ell-1)}\frac{[n\ell]_{\xi}}{[\ell]_{\xi}} K_i^{\ell} k_{i,n\ell} 
= k_{i,n\ell} \Big( k_{i,\ell} 
- n K_i^{\ell} \Big).
\end{align*}
Hence $k_{i, (n+1)\ell} = \frac{1}{n+1} k_{i,n\ell} \Big( k_{i,\ell} - n K_i^{\ell} \Big)$, so we obtain \eqref{eq:K0nl-formula} recursively.

Finally we take $t=m\ell$, $t'\in\I_{\ell - 1}$ in \eqref{eq:g8}. Using \eqref{eq:binom-evaluation-2},
\begin{align*}
k_{i,m\ell+t'} & = 
\sum_{0 \le j \le t'} (-1)^j 
{\pinom{m\ell+j-1}{j}}_{\xi} K_i^j k_{i,m\ell}k_{i,t'-j} = k_{i,m\ell}k_{i,t'}.
\end{align*}
Hence the claim follows from Proposition \ref{prop:V0-prop-lusztig} \ref{item:V0-prop-lusztig-1}.
\epf

Let $\Gamma = (\Z/{2\ell})^{\I}$, with $g_i\in\Gamma$ being 
generators of the corresponding copies of the cyclic group $\Z/{2\ell}$. Let $\hgo$ be the abelian Lie algebra with basis  $(t_i)_{i\in\I}$, so that $U(\hgo) \simeq \cB[t_i: i\in\I]$.
 
\begin{theorem}\label{th:V0-isom} \cite[Theorem 4.1]{Lentner2}
The assignment
\begin{align}\label{eq:V0-isom}
\Psi(g_i)&=K_i, & \Psi(t_i)&=h_{i,\ell}, & i&\in \I,
\end{align}
determines an isomorphism of Hopf algebras  $\Psi: \cB\Gamma \otimes U(\hgo) \to V^0_{\cB}$.
\end{theorem}

We present a different proof involving the polynomials $p_{i,t}$.

\pf That \eqref{eq:V0-isom} defines an algebra map follows by \eqref{eq:g6} and \eqref{eq:K2l=1}; 
that is a surjective Hopf algebra map, by \eqref{eq:comult-V-3}, 
Definition \ref{def:hn-skew-primitive} and Lemma \ref{lemma:V0-generators}. 
It remains to prove that $\Psi$ is injective.
By \cite[5.3.1]{Mo-libro} it reduces to prove that 
$\Psi$ is injective on the first term of the coradical filtration, i.e. that
the set $\{K_i^p h_{r,\ell}^j: i,r\in\I_\theta, p\in\I_{0,2\ell-1}, j\in\I_{0,1} \}$ is linearly independent. 
By the assumption $\xi^2 \neq 1$, we have $p_{\ell,\ell}(\xi)=n$, so \eqref{eq:defn-hn} implies that
\begin{align}\label{eq:hl=kl+trash}
h_{i,\ell} &\in k_{i,\ell} K_i^{\ell} + \cB \langle K_i\rangle,
\end{align}
see the line before \eqref{eq:K0nl-formula}.
We need then to prove that the set 
\begin{align*}
\{K_i^p k_{r,\ell}^j: i,r\in\I_\theta, p\in\I_{0,2\ell-1}, j\in\I_{0,1} \}
\end{align*}
is linearly independent. 
Indeed, suppose that 
\begin{align}\label{eq:V0-isom-inj}
0=\sum_{i,p} e_{i,p} K_i^p +b_{i,p} K_i^p k_{i,\ell},
\end{align} 
where $e_{i,p},b_{i,p}\in\cB$.
Fix $i\in \I$. The $\cA'$-algebra maps $\Xi_{i,j}: \cA' [\Z] \to \cA'$ as in \eqref{eq:defn-xi-character} 
satisfy $\Xi_{i,j}(V^0) \subseteq \cA$ by \eqref{eq:maps-preserve}. We restrict to $\cA$-algebra maps $\Xi_{i,j}: V^0 \to \cA$ 
and tensorize  to get $\cB$-algebra maps
$\Xi_{i,j}: V^0_{\cB} \to \cB$ such that 
\begin{align*}
\Xi_{i,j}(K_i) &= \xi^j,& \Xi_{i,j}(k_{i,\ell}) &= {\pinom{j}{\ell}}_{\xi},& \Xi_{i,j}(K_r) &= 1,& 
\Xi_{i,j}(k_{r,\ell})&=0 \text{ if } r\neq i.
\end{align*} Applying $\Xi_{i,j}$ to \eqref{eq:V0-isom-inj}, we get
\begin{align} \label{eq:linear-comb}
0&=\sum_{p\in\I_{0,2\ell-1}} e_{i,p} \xi^{pj}, &
0&\le j<\ell;
\\ \label{eq:linear-comb-2}
0&=\sum_{p\in\I_{0,2\ell-1}} e_{i,p} \xi^{pj}+b_{i,p} \xi^{pj}, &
\ell & \le j<2\ell.
\end{align}

If $\ell'$ is even, then $\ell'=2\ell$ and from \eqref{eq:linear-comb} we deduce that $e_{i,p}=0$ for all $p\in\I_{0,2\ell-1}$. Hence
$0=\sum_{p\in\I_{0,2\ell-1}} b_{i,p} \xi^{pj}$ for all $0\le j<\ell$ by \eqref{eq:linear-comb-2}, and the same argument shows that $b_{i,p}=0$ for all $p\in\I_{0,2\ell-1}$.

If $\ell'=\ell$ is odd, then $e_{i,p}+e_{i,p+\ell}\overset{\star}{=} 0$ for all $p\in\I_{0,\ell-1}$ by \eqref{eq:linear-comb}.
Similarly as above, we consider the algebra maps $\widetilde{\Xi}_{i,j}: \cA' [\Z^j] \to \cA'$ such that $K_i\mapsto -v^j$ and $K_r\mapsto 1$ for $r\neq i$; 
we get algebra maps $\widetilde{\Xi}_{i,j}: V^0_{\cB} \to \cB$ such that
\begin{align*}
\widetilde{\Xi}_{i,j}(K_i)&=-\xi^j, & \widetilde{\Xi}_{i,j}(k_{i,\ell})&=-{\pinom{j}{\ell}}_{\xi}, & 
\Xi_{i,j}(K_r)&=1, & \Xi_{i,j}(k_{r,\ell})&=0, & r&\neq i. 
\end{align*}
Applying $\widetilde{\Xi}_{i,j}$ to the previous equality \eqref{eq:V0-isom-inj}, we see that
\begin{align*}
0&=\sum_{p\in\I_{0,\ell-1}} (-1)^{i} (e_{i,p}-e_{i,p+\ell}) \xi^{pj}, &
0&\le j<\ell.
\end{align*}
Hence $e_{i,p}-e_{i,p+\ell}=0$ for all $p\in\I_{0,\ell-1}$, so $e_{i,p}=0$ for all $p\in \I_{0,2\ell-1}$. Analogously $b_{i,p}=0$ for all $p\in \I_{0,2\ell-1}$.
\epf

\section{The algebra $V$, simply-laced diagram}\label{sec:V-slaced}
\subsection{Definitions and first properties} 
As in \cite{L-fdHa-JAMS}, we fix a finite Cartan matrix $A=(a_{ij})_{i,j\in\I}$ whose Dynkin diagram is connected and simply-laced, that is, of type A, D or E.

\medbreak
Following \cite[2.3, pp. 268 ff]{L-fdHa-JAMS} we consider the $\cA$-algebra
$V$ presented by generators \eqref{eq:V0-gens}, $E_i^{(N)}$, $F_i^{(N)}$, 
$i\in \I$, $N \in \N_0$
with relations  \eqref{eq:g5}, \dots  \eqref{eq:g10}, together with the following, tagged again as in \emph{loc. cit.},
\begin{align}
\tag{d1} E_i^{(N)}E_i^{(M)} &= {\pinom{N+M}{M}}_{v} E_i^{(N + M)},  & E_i^{(0)} &= 1; \label{eq:d1}
\\
\tag{f1} F_i^{(N)}F_i^{(M)} &= {\pinom{N+M}{M}}_{v} F_i^{(N + M)}, & F_i^{(0)} &= 1; \label{eq:f1}
\end{align}
if $i\neq j\in\I$, $a_{ij}=0$:
\begin{align}
\tag{d2} E_i^{(N)}E_j^{(M)} &= E_j^{(M)}E_i^{(N)}, \label{eq:d2}
\\
\tag{f2} F_i^{(N)}F_j^{(M)} &= F_j^{(M)}E_i^{(N)}, \label{eq:f2}
\end{align}
if $i\neq j\in\I$, $a_{ij}=-1$, $i<j$:
\begin{align}
\tag{d3} E_i^{(N)}E_j^{(M)} &= \sum_{t=0}^{\min \{M,N\}} v^{t+(N-t)(M-t)}E_j^{(M-t)}E_{ij}^{(t)}E_i^{(N-t)}, \label{eq:d3}
\\
\tag{d4} v^{NM}E_i^{(N)}E_{ij}^{(M)} &= E_{ij}^{(M)}E_i^{(N)}, \label{eq:d4}
\\
\tag{d5} v^{NM}E_{ij}^{(M)}E_j^{(N)} &= E_j^{(N)}E_{ij}^{(M)}, \label{eq:d5}
\\
\tag{f3} F_i^{(N)}F_j^{(M)} &= \sum_{t=0}^{\min \{M,N\}} v^{-t-(N-t)(M-t)} F_j^{(M-t)}F_{ij}^{(t)}F_i^{(N-t)}, \label{eq:f3}
\\
\tag{f4} v^{NM}F_i^{(N)}F_{ij}^{(M)} &= F_{ij}^{(M)}F_i^{(N)}, \label{eq:f4}
\\
\tag{f5} v^{NM}F_{ij}^{(M)}F_j^{(N)} &= F_j^{(N)}F_{ij}^{(M)}, \label{eq:f5}
\end{align}
where $E_{ij}^{(N)}= \sum\limits_{k=0}^{N} (-1)^{N-k}v^{-k}E_i^{(k)} E_j^{(N)} E_i^{(N-k)}$ (cf. \cite[Lemma 2.5 (d)]{L-fdHa-JAMS}), $F_{ij}^{(N)}= \sum\limits_{k=0}^{N} (-1)^{N-k}v^{-k}F_i^{(k)} F_j^{(N)} F_i^{(N-k)}$;
\begin{align}
\tag{h1} E_i^{(N)}F_j^{(M)} &= F_j^{(M)} E_i^{(N)}, \qquad\qquad i\neq j, \label{eq:h1}
\\ \tag{h2} E_i^{(N)}F_i^{(M)} &= \sum_{0 \le t \le \min\{N,M\}} F_i^{(M - t)} \pinom{K_i; 2t -N - M}{t} E_i^{(N - t)}, \label{eq:h2}
\\ \tag{h3} K_i^{\pm 1}E_j^{(N)}  &= v^{\pm Na_{ij}}   E_j^{(N)} K_i^{\pm 1}, \label{eq:h3}
\\ \tag{h4} K_i^{\pm 1}F_j^{(N)}  &= v^{\mp Na_{ij}}   F_j^{(N)} K_i^{\pm 1}, \label{eq:h4}
\\ \tag{h5} \pinom{K_i;c}{t} E_j^{(N)}  &= E_j^{(N)} \pinom{K_i;c + Na_{ij}}{t}, \label{eq:h5}
\\ \tag{h6} \pinom{K_i;c}{t} F_j^{(N)}  &= F_j^{(N)} \pinom{K_i;c - Na_{ij}}{t}. \label{eq:h6}
\end{align}
Let $V^+$, respectively $V^-$, be the subalgebra of $V$
generated by $E_i^{(N)}$, respectively $F_i^{(N)}$, $i\in\I$, $N \in \N_0$.
Let 
\begin{align}\label{eq:pinom-Ki-1}
\pinom{K_i^{-1};c}{t} &=\Ss\left(\pinom{K_i;c}{t}\right). 
\end{align}
The following formula is analogous to \eqref{eq:h2}, cf. \cite[Corollary 3.19]{L-libro}: 
\begin{align}
F_i^{(N)}E_i^{(M)} &= \sum_{0 \le t \le \min\{N,M\}} E_i^{(M - t)} \pinom{K_i^{-1}; 2t -N - M}{t} F_i^{(N - t)}. \label{eq:h2-v2}
\end{align}

By \cite[Proposition 4.8, p. 287]{L-fdHa-JAMS}, we know that $V$ has a unique Hopf algebra structure 
determined by \eqref{eq:comult-V-3} and
\begin{align}\label{eq:comult-V}
&\begin{aligned}
\Delta (E_i^{(N)}) &=  \sum_{0 \le b \le N}   v^{b(N - b)} E_i^{(N - b)}K_i^b \otimes  E_i^{(b)}, 
\\ \Delta (F_i^{(N)}) &= \sum_{0 \le a \le N}   v^{-a(N - a)} F_i^{(a)} \otimes K_i^{-a} F_i^{(N - a)},
\end{aligned}
& i&\in \I, \, N \in \N_0.
\end{align}

\subsection{Specializations of $V$}\label{subsec:specializations-slaced}

We define next \begin{align*}
V^+_{\cB} &= V^+ \otimes_\cA \cB,& V^-_{\cB} &= V^- \otimes_\cA \cB,& V_{\cB} &= V \otimes_\cA \cB.
\end{align*}
By \cite[Proposition 3.2 (b)]{L-contemporary}, 
$V_{\cB}^+$ is generated by $E_i$ and $E_i^{(\ell)}$, $i\in \I$; $V_{\cB}^-$ is generated by $F_i$ and $F_i^{(\ell)}$, $i\in \I$.
From now on, we abbreviate
\begin{align*}
k_{i,N} &= \pinom{K_i;0}{N},& N \in \N_0,& & E_i &= E_i^{(1)}, 
& F_i &= F_i^{(1)}.
\end{align*}

\begin{lemma} Let $i, j\in\I$. We have in $V_{\cB}$:
\begin{align}
\label{eq:kiell-Ei}
k_{i,\ell} E_i &= E_i\big(k_{i,\ell} + \xi^{-2} [2]_{\xi} K_i^{-1} k_{i,\ell - 1} + \xi^{-4} K_i^{-2} k_{i,\ell - 2}  \big),
\\
\label{eq:kiell-Fi}
k_{i,\ell} F_i &= F_i \Big( k_{i,\ell} + \sum_{j \in \I_{\ell-2}} (-1)^j \xi^{-2j} [j+1]_{\xi} K_i^j k_{i,\ell-j} - K_i^{\ell}\Big),
\\
\label{eq:kiell-Ej-aij=0}
k_{i,\ell} E_j &= E_j k_{i,\ell}   , \qquad i\neq j, \, a_{ij}=0, 
\\
\label{eq:kiell-Fj-aij=0}
k_{i,\ell} F_j &=F_j k_{i,\ell}, \qquad i\neq j, \, a_{ij}=0,
\\
\label{eq:kiell-Ej-aij=-1}
k_{i,\ell} E_j &= \xi^{\ell} E_j\Big( \sum_{s=0}^{\ell}
(-\xi)^{-s} K_i^{s}k_{i,\ell-s} \Big),  \qquad i\neq j, \, a_{ij}=-1,
\\
\label{eq:kiell-Fj-aij=-1}
k_{i,\ell} F_j &= \xi^{\ell} F_j\left( k_{i,\ell}
+\xi^{-1} K_i^{-1}k_{i,\ell-1}\right),  \qquad i\neq j, \, a_{ij}=-1,
\\
\label{eq:kiell-Ejell}
k_{i,\ell} E_j^{(\ell)} &= E_j^{(\ell)} \left(k_{i,\ell} + a_{ij} K_i^{\ell} \right),
\\
\label{eq:kiell-Fjell}
k_{i,\ell} F_j^{(\ell)} &= F_j^{(\ell)} \left(k_{i,\ell} - a_{ij} K_i^{\ell} \right).
\end{align}
\end{lemma}
\pf
We consider first the case $j=i$. We take $t=\ell$, $c=0$ and $N=1$ in \eqref{eq:h6} and use \eqref{eq:g9} to obtain \eqref{eq:kiell-Fi}:
\begin{align*}
k_{i,\ell} F_i &= F_i \pinom{K_i;-2}{\ell} = F_i \Big( \sum_{0 \le j \le \ell} (-1)^j \xi^{2(\ell-j)} 
{\pinom{j+1}{j}}_{\xi} K_i^j k_{i,\ell-j} \Big)
\\
&= F_i \Big( k_{i,\ell} + \sum_{j \in \I_{\ell-2}} (-1)^j \xi^{-2j} [j+1]_{\xi} K_i^j k_{i,\ell-j} +(-\xi)^{\ell} K_i^{\ell}\Big).
\end{align*}
For \eqref{eq:kiell-Fjell}, we take $t=\ell$, $c=0$ and $N=\ell$ in \eqref{eq:h6} and use \eqref{eq:g9}:
\begin{align*}
k_{i,\ell} F_i^{(\ell)} &= F_i^{(\ell)} \pinom{K_i;-2\ell}{\ell} = F_i^{(\ell)} \Big( \sum_{0 \le j \le \ell} (-1)^j \xi^{2\ell(\ell-j)} {\pinom{2\ell+j-1}{j}}_{\xi} K_i^j k_{i,\ell-j} \Big)
\\
&= F_i^{(\ell)} \left(k_{i,\ell}+(-1)^{\ell}\xi^{2\ell(\ell-1)} 2 K_i^{\ell} \right)
= F_i^{(\ell)} \left(k_{i,\ell} - 2 K_i^{\ell} \right).
\end{align*}

Now we take $j\neq i$ and $a_{ij}=0$. From \eqref{eq:h6},
$k_{i,\ell}F_j^{(N)}=F_j^{(N)}k_{i,\ell}$ for all $N\in\N$, hence we obtain \eqref{eq:kiell-Fj-aij=0} when $N=1$, and \eqref{eq:kiell-Fjell} when $N=\ell$.

Next we take $j\neq i$ and $a_{ij}=-1$. From \eqref{eq:h6} and \eqref{eq:g10} we derive \eqref{eq:kiell-Fj-aij=-1} when $N=1$, and \eqref{eq:kiell-Fjell} when $N=\ell$:
\begin{align*}
k_{i,\ell}F_j&= F_j \pinom{K_i;1}{\ell} = F_j\left( \sum_{s=0}^{\ell}
\xi^{\ell-s} {\pinom{1}{s}}_{\xi} K_i^{-s}\pinom{K_i; 0}{\ell-s} \right)
\\
& = \xi^{\ell}F_j\left( k_{i,\ell}
+\xi^{-1} K_i^{-1}k_{i,\ell-1}\right),
\\
k_{i,\ell}F_j^{(\ell)}&= F_j^{(\ell)} \pinom{K_i;\ell}{\ell} =
F_j\left( \sum_{s=0}^{\ell}
\xi^{\ell(\ell-s)} {\pinom{\ell}{s}}_{\xi} K_i^{-s}\pinom{K_i; 0}{\ell-s} \right)
 = F_j\left( k_{i,\ell} + K_i^{-\ell}\right).
\end{align*}

Finally we get \eqref{eq:kiell-Ei}, \eqref{eq:kiell-Ej-aij=0},
\eqref{eq:kiell-Ej-aij=-1} and \eqref{eq:kiell-Ejell} similarly but from \eqref{eq:h5}.
\epf

Now we compute relations involving $h_{i,\ell}$.
The formulas \eqref{eq:hellE} and \eqref{eq:hellF} appear in \cite[Theorem 3.1]{Lentner}.

\begin{lemma}\label{lem:action-primitive} Let $i,j\in\I$. We have in $V_{\cB}$:
\begin{align}
\label{eq:hellE}
h_{i,\ell} E_j &= E_jh_{i,\ell} + a_{ij}\ell^{-1} E_j,
\\
\label{eq:hellEell}
h_{i,\ell} E_j^{(\ell)} &= E_j^{(\ell)}h_{i,\ell} + a_{ij} E_j^{(\ell)},
\\
\label{eq:hellF}
h_{i,\ell} F_j &= F_j h_{i,\ell}- a_{ij}\ell^{-1},
\\
\label{eq:hellFell}
h_{i,\ell} F_j^{(\ell)} &= F_j^{(\ell)}h_{i,\ell} - a_{ij} F_j^{(\ell)}.
\end{align}
\end{lemma}
\pf
By \eqref{eq:defn-hn}, there exist $\mathtt{b}_t \in \cB$ such that
\begin{align}\label{eq:hl=kl+grouplike}
h_{i,\ell}&=
\Big(\sum_{s \in \I_\ell} \frac{p_{\ell,s}}{\ell} \, k_{i,s} K_i^{\phi_\ell(s)} \Big) K_i^{\ell}
=
k_{i,\ell} K_i^{\ell} + \sum_{t=0}^{\ell-1} \mathtt{b}_t K_i^{2t}.
\end{align}
Indeed, for each $s\in\I_{\ell-1}$,
\begin{align*}
k_{i,s} &= \frac{1}{[s]_{\xi}^!} \prod_{j=0}^{s-1} \frac{\xi^{-j}K_i-\xi^jK_i^{-1}}{\xi-\xi^{-1}} \in \sum_{p=0}^s \cB K_i^{2p-s}.  
\end{align*}

Using \eqref{eq:h3}, \eqref{eq:kiell-Ejell} and \eqref{eq:K2l=1},
\begin{align*}
h_{i,\ell}& E_j^{(\ell)} = (k_{i,\ell} K_i^{\ell} + \sum_{t=0}^{\ell-1} \mathtt{b}_t K_i^{2t}) E_j^{(\ell)}
\\
= & \xi^{a_{ij}\ell^2}E_j^{(\ell)} \left(k_{i,\ell} + a_{ij} K_i^{\ell} \right) K_i^{\ell} +  \sum_{t=0}^{\ell-1} \mathtt{b}_t \xi^{2t a_{ij}\ell} E_j^{(\ell)} K_i^{t}
= E_j^{(\ell)}h_{i,\ell} + a_{ij} E_j^{(\ell)}.
\end{align*}

The proof of \eqref{eq:hellFell} is similar.
Next we check \eqref{eq:hellE}. By a direct computation,
\begin{align*}
\Delta ([h_{i,\ell}, E_j]) &= [h_{i,\ell}, E_j] \otimes 1 + K_j \ot [h_{i,\ell}, E_j].
\end{align*}
Thus $[h_{i,\ell}, E_j]$ is $(1,K_j)$-primitive and belongs to the subalgebra generated by $K_j$, $h_{i,\ell}$ and $E_j$, so $[h_{i,\ell}, E_j] = c_{ij} E_j+ d_{ij} (1-K_j)$ for some 
$c_{ij},d_{ij}\in\cB$. As $E_j^\ell=0$, we have that
\begin{align*}
0 & = [h_{i,\ell}, E_j^\ell] = \sum_{k\in\I_{\ell}} E_j^{k-1} [h_{i,\ell}, E_j] E_j^{\ell-k}
\\
& = \sum_{k\in\I_{\ell}} E_j^{k-1} (c_{ij}E_j+ d_{ij} (1-K_j)) E_j^{\ell-k}
= d_{ij} \sum_{k\in\I_{\ell}} E_j^{k-1}  (1-K_j) E_j^{\ell-k}
\\ 
& = d_{ij} \Big( \ell E_j^{\ell-1}  - \big(\sum_{k\in\I_{\ell}}\xi^{-2k}\big) E_j^{\ell-1} K_j \Big)
= d_{ij} \ell E_j^{\ell-1}.
\end{align*}
Hence $d_{ij}=0$.
Analogously $[h_{i,\ell}, F_j] = c_{ij}'F_j$ for some $c_{ij}'\in\cB$. Now
\begin{align*}
0 &=\left[ h_{i,\ell}, K_j-K_j^{-1} \right] 
= (\xi-\xi^{-1}) \left[ h_{i,\ell}, [E_j,F_j] \right]
= (c_{ij}+c_{ij}') \left(K_j-K_j^{-1}\right),
\end{align*}
so $c_{ij}'=-c_{ij}$. 
We consider three cases:
\begin{enumerate}[leftmargin=*]
\item $j\neq i$, $a_{ij}=0$. Then $[h_{i,\ell},E_j]=0$ by \eqref{eq:kiell-Ej-aij=0} and \eqref{eq:h3}. Thus $c_{ij}=0$.
\item $j\neq i$, $a_{ij}=-1$. From \eqref{eq:kiell-Fj-aij=-1}, \eqref{eq:hl=kl+grouplike} and \eqref{eq:h4}:
\begin{align*}
-c_{ij}F_j &= [h_{i,\ell},F_j] = [k_{i,\ell}K_i^{\ell},F_j] + \sum_{t=0}^{\ell-1} \mathtt{b}_t [K_i^{2t},F_j]
\\
&=\xi^{-1}F_jK_{i}^{\ell-1}k_{i,\ell-1}+\sum_{t=0}^{\ell-1}\mathtt{b}_t(\xi^{2t}-1)F_jK_i^{2t}.\\
\end{align*}
As the set $\{F_j K_i^{2t}: t\in\I_{0,\ell-1} \}$ is linearly independent, 
$$c_{ij}=-\xi^{-1}b_{0}$$
where $b_t\in\cB$ denote the elements satisfying
\begin{align*}
k_{i,\ell-1} &= \frac{1}{[\ell-1]_{\xi}^!} \prod_{j=0}^{\ell-2} \frac{\xi^{-j}K_i-\xi^jK_i^{-1}}{\xi-\xi^{-1}}=\sum_{t=0}^{\ell-1} b_t K_i^{2t-\ell+1}.
\end{align*}
Since 
\begin{align*}
b_0= \frac{1}{[\ell-1]_{\xi}^!}  \frac{\prod_{j=0}^{\ell-2}-\xi^j}{(\xi-\xi^{-1})^{\ell-1}}
= \frac{(-1)^{\ell-1}\xi^{\binom{\ell-1}{2}}}{\xi^{\binom{\ell}{2}}\prod_{j\in\I_{\ell-1}}1-\xi^{-2j}}
=\frac{\xi}{\ell},
\end{align*}
we have that 
$c_{ij}=-\ell^{-1}$.

\item $j=i$. The proof is analogous to the previous case, using \eqref{eq:kiell-Ei}.
\end{enumerate}

\medspace

Hence $c_{ij}=\frac{a_{ij}}{\ell}$ in all the cases, so \eqref{eq:hellE} and \eqref{eq:hellF} follow.
\epf

\bigbreak

\subsection{The Hopf algebra structure of $V_{\cB}$}
Recall that by Theorem \ref{th:V0-isom}, $V_{\cB}^0 = \cB[K_i, h_{i,\ell}: i\in \I] \simeq \cB\Gamma \otimes U(\hgo)$.

\begin{remark}\label{rem:counit-Ki-c-t}
The counit on the elements $\displaystyle\pinom{K_i^{\pm1};c}{t}$ takes the following values.
\begin{align}\label{eq:counit Ki c t}
	\varepsilon\left(\pinom{K_i^{\pm1};c}{t}\right)
	=\begin{cases} 
	1&\mbox{if $c=t=0$,}\\
	0&\mbox{if $c=0$ and $t\neq 0$,}\\
	\pinom{c}{t}_\xi&\mbox{if $c>0$,}\\
	(-1)^t\pinom{-c+t-1}{t}_\xi&\mbox{if $c<0$.}\\
	\end{cases}
	\end{align}
	In fact, we first note that $\displaystyle \varepsilon\left(\pinom{K_i^{-1};c}{t}\right)=
	\varepsilon\Ss\left(\pinom{K_i;c}{t}\right)=
	\varepsilon\left(\pinom{K_i;c}{t}\right)$
	by \eqref{eq:pinom-Ki-1}. The formula for $c=0$ holds by \eqref{eq:kt}. Then, for $c>0$, we use \eqref{eq:g10}:
	$$
	\varepsilon\left(\pinom{K_i; c}{t}\right)= 
	\sum_{0 \le j \le t}  v^{c(t-j)} {\pinom{c}{j}}_{v} \varepsilon( K_i^{-j}) \varepsilon\left(\pinom{K_i; 0}{t - j}\right)=\pinom{c}{t}_v.
	$$
	While for $c<0$, we use \eqref{eq:g9}:
	\begin{align*}
	\varepsilon\left(\pinom{K_i; c}{t}\right)&= 
	\sum_{0 \le j \le t} (-1)^j v^{-c(t-j)} {\pinom{-c+j-1}{j}}_{v} \varepsilon( K_i^j) \varepsilon\left(\pinom{K_i; 0}{t - j}\right)
	\\
	&=(-1)^t{\pinom{-c+t-1}{t}}_{v} 
	\end{align*}
	\end{remark}

\begin{theorem}\label{th:structure-simply-laced}
The Hopf algebra
$V_{\cB}$ has a triangular decomposition given by a TD-datum $( V_{\cB}^+, V_{\cB}^-, \rightharpoonup, \leftharpoonup, \sharp)$ over $V_{\cB}^0$. 
The left action $\rightharpoonup$ of $V_{\cB}^-$ on $V_{\cB}^+$, the right
action $\leftharpoonup$ of $V_{\cB}^+$ on $V_{\cB}^-$ and the map $\sharp:V_{\cB}^- \otimes V_{\cB}^+ \to V_{\cB}^0$
are determined as follows:
\begin{align}\label{eq:action-pairing-V-V+}
\begin{aligned}
F_i^{(N)}\rightharpoonup E_j^{(M)} &=  \delta_{ij}
(-1)^N\pinom{M-1}{N}_\xi E_i^{(M-N)}
\\
F_i^{(N)}\leftharpoonup E_j^{(M)} &= \delta_{ij} (-1)^{M} {\pinom{N-1}{M}}_{\xi}  F_i^{(N - M)}
\\
F_i^{(N)} \sharp E_j^{(M)} &= \delta_{M,N}\delta_{ij}  \pinom{K_i^{-1};0}{N},
\end{aligned}\end{align}
cf. \eqref{eq:pinom-Ki-1}, where $E_i^{(n)}=0=F_i^{(n)}$ if $n<0$. 
\end{theorem}

\pf For  the first claim, we just need to verify that the conditions of Proposition \ref{prop:sommerh} \ref{item:prop-sommerh-2} hold. 

\smallbreak
Let $V_{\cB}^{\geq 0} := V_{\cB}^+V_{\cB}^0$ and $V_{\cB}^{\leq 0} := V_{\cB}^0V_{\cB}^-$; these are Hopf subalgebras of $V_{\cB}$ by definition. 
It is easy to see 
that the inclusions $V_{\cB}^0 \hookrightarrow V_{\cB}^{\geq 0}$ and $V_{\cB}^0 \hookrightarrow V_{\cB}^{\leq 0}$ 
admit  Hopf algebra sections $\pi^+$ and $\pi^-$ respectively
and that
\begin{align*}
V_{\cB}^+ &= \left(V_{\cB}^{\geq 0}\right)^{\co \pi^+}, & V_{\cB}^- &= {}^{\co \pi^-}\hspace{-5pt}\left(V_{\cB}^{\geq 0}\right).
\end{align*}
Thus $V_{\cB}^+$ is a Hopf algebra in $\ydV$ and $V_{\cB}^{\geq 0} \simeq V_{\cB}^+ \# V_{\cB}^0$, respectively
$V_{\cB}^-$ is a Hopf algebra in $\Vyd$ and $V_{\cB}^{\leq 0} \simeq V_{\cB}^0 \# V_{\cB}^-$.
Also, by \cite[Theorem 4.5 (a)]{L-fdHa-JAMS}, the multiplication induces a linear isomorphism $V_{\cB}^+ \otimes V_{\cB}^0 \otimes V_{\cB}^- \simeq V_{\cB}$.
Thus we may apply  Proposition \ref{prop:sommerh} \ref{item:prop-sommerh-2}. 

The verification of \eqref{eq:action-pairing-V-V+} is direct using the formulas in the proof of \cite[Theorem 3.5]{S} and the natural projections
$\varpi^{\star}: V_{\cB} \to V_{\cB}^{\star}$, for $\star \in \{+, 0, -\}$. In fact, $F_i^{(N)}\rightharpoonup E_j^{(M)} =
\varpi^{+}(F_i^{(N)}E_j^{(M)})$. If $i\neq j$, this zero by \eqref{eq:h1}. Otherwise, we use \eqref{eq:h2-v2}:
\begin{align*}
F_i^{(N)}\rightharpoonup E_i^{(M)} &=
\varpi^{+}(F_i^{(N)}E_i^{(M)})\\
&= \sum_{0 \le t \le \min\{N,M\}} E_i^{(M - t)} \varepsilon\left(\pinom{K_i^{-1}; 2t -N - M}{t}\right) \varepsilon\left(F_i^{(N - t)}\right)
\end{align*}
which is zero for $N\geq M$ by Remark \ref{rem:counit-Ki-c-t}. If $N<M$, then
\begin{align*}
F_i^{(N)}\rightharpoonup E_i^{(M)} &=E_i^{(M-N)}\varepsilon\left(\pinom{K_i^{-1}; N - M}{N}\right)=(-1)^N\pinom{M-1}{N}E_i^{(M-N)}.
\end{align*}

We can verify the formulas for $\leftharpoonup$ and $\sharp$ in a similar way.

We next show by induction that \eqref{eq:action-pairing-V-V+} completely determines $\rightharpoonup$, $\leftharpoonup$ and $\sharp$. We will use that the comultiplication of $V_{\cB}^{\pm}$ in the respective Yetter-Drinfeld category is given by
\begin{align*}
&\begin{aligned}
\Delta (E_i^{(N)}) &=  \sum_{0 \le b \le N}   v^{b(N - b)} E_i^{(N - b)} \otimes  E_i^{(b)}, 
\\ \Delta (F_i^{(N)}) &= \sum_{0 \le a \le N}   v^{-a(N - a)} F_i^{(a)} \otimes F_i^{(N - a)},
\end{aligned}
& i&\in \I, \, N \in \N_0.
\end{align*}
This follows from \eqref{eq:comult-V}.

Let $E=E_{j_1}^{(M_1)}\cdots E_{j_r}^{(M_r)}$ and $F=F_{i_1}^{(N_1)}\cdots F_{i_s}^{(N_s)}$. First, we assume that 
\begin{align*}
F_i^{(N)}\rightharpoonup E,\quad F\leftharpoonup E_j^{(M)},\quad F_i^{(N)}\sharp E\quad\mbox{ and }\quad F\sharp E_j^{(M)}
\end{align*}
are determined by \eqref{eq:action-pairing-V-V+} for all $r,s\leq n$ and prove the same claim for $s=r=n+1$. By \eqref{eq:Compatibility harpoonup}, we have that
\begin{align}
\notag
F_i^{(N)}\rightharpoonup\left(E_j^{(M)}E\right)&=\left({\left(F_i^{(N)}\right)\^{1}}\_{0} \rightharpoonup \left(E_j^{(M)}\right)\^{1}\right)\times
\\ \label{eq:F rightharpoonup EE}
\times
\left(
{\left(F_i^{(N)}\right)\^{1}}\_{1}
\right.&\left.
\left(\left(F_i^{(N)}\right)\^{2} \sharp \left(E_j^{(M)}\right)\^{2}\right) {\left(E_j^{(M)}\right)\^{3}}\_{-1}
\right)\actizq
\\
\notag& \left(\left(\left(F_i^{(N)}\right)\^{3} \leftharpoonup {\left(E_j^{(M)}\right)\^{3}}\_{0}\right) \rightharpoonup E\right).
\end{align}
Hence, \eqref{eq:F rightharpoonup EE} is determined by \eqref{eq:action-pairing-V-V+} because of the inductive hypothesis. The same holds for $F_i^{(N)}\sharp\left(E_j^{(M)}E\right)$ since 
\begin{align*}
F_i^{(N)}\sharp\left(E_j^{(M)}E\right)&= \left(\left(F_i^{(N)}\right)\^{1} \sharp \left(E_j^{(M)}\right)\^{1}\right)\times\\ &{\left(E_j^{(M)}\right)\^{2}}\_{-1}\left(\left(\left(F_i^{(N)}\right)\^{2} \leftharpoonup {\left(E_j^{(M)}\right)\^{2}}\_{0}\right) \sharp E\right)
\end{align*}
by \eqref{eq:Compatibility sharp}. A similar argument works for $\left(FF_i^{(N)}\right)\leftharpoonup E_j^{(M)}$ and $\left(FF_i^{(N)}\right)\sharp E_j^{(M)}$.

Second, we prove that 
\begin{align*}
F\rightharpoonup E,\quad F\leftharpoonup E\quad\mbox{ and }\quad F\sharp  
E
\end{align*}
are determined by \eqref{eq:action-pairing-V-V+} for all $r,s\geq0$. This is true for $F\rightharpoonup E$ and $F\leftharpoonup E$ because $\rightharpoonup$ and $\leftharpoonup$ are actions. For the others we proceed again by induction on $r$ (or on $s$) using \eqref{eq:Compatibility sharp}; notice that the initial inductive step $r=1$ was proved above.
\epf

\begin{remark} Here are some particular instances of the first line in \eqref{eq:action-pairing-V-V+}:
\begin{align*}
F_i \rightharpoonup E_j^{(M)} &= 
F_i^{(\ell)} \rightharpoonup E_j^{(M)}=0,
\qquad\qquad\qquad \text{if }i\neq j,
\\
F_i \rightharpoonup E_i^{(M)} &=
(-1)^{M-1} [M-1]_{\xi} E_i^{(M - 1)},
\\
F_i^{(\ell)} \rightharpoonup E_i^{(M)} &= 0 \qquad\qquad\qquad \text{if }\ell \text{ does not divide } M,
\\
F_i^{(\ell)} \rightharpoonup E_i^{(\ell n)} &= (-1)^{n-1} (n-1) E_i^{(\ell n - \ell)}.
\end{align*}
\end{remark}

\begin{remark} The structure of $V_{\cB}^+$ as an object in $\ydV$ is  as follows:
	the (left) action of $V_{\cB}^0$ on $V_{\cB}^+$ is given by \eqref{eq:h3}, \eqref{eq:hellE} and \eqref{eq:hellEell},
	while the coaction $\lambda: V_{\cB}^+ \to V_{\cB}^0 \otimes V_{\cB}^+$ is determined by
	\begin{align*}
	\lambda(E_i^{(N)}) &= K_i^{N}\ot E_i^{(N)},& i\in \I, \ N &\in \N.
	\end{align*}
	Analogously, the structure of $V_{\cB}^-$ as an object in $\Vyd$ is as follows:
	the (right) action of $V_{\cB}^0$ on $V_{\cB}^-$ is given by \eqref{eq:h4}, \eqref{eq:hellF} and \eqref{eq:hellFell};
	meanwhile the coaction $\rho: V_{\cB}^- \to V_{\cB}^- \otimes V_{\cB}^0$  is determined by
	\begin{align*}
	\rho(F_i^{(N)}) &= F_i^{(N)}\ot K_i^{-N},& i\in \I, \ N &\in \N.
	\end{align*}
\end{remark}

\subsection{The multiply-laced diagrams}\label{subsec:V-mlaced}
The arguments above can be extended to the diagrams of types B, C, F, G. We just discuss the torus part here.

Let $\bd = (d_i)_{i\in \I} \in \N^{\I}$.
Following \cite[6.4]{L-qgps-at-roots} we consider the $\cA$-algebra $\V^0$ that is a (multiply-laced!) variation of 
the $V^0$ studied so far. For the agility of the exposition we do not stress $\bd$ in the notation.
This $\V^0$ is presented by the generators analogous to those \eqref{eq:V0-gens} of $V^0$:
\begin{align}\label{eq:V0-gens-multlaced}
&\K_i, & &\K_i^{-1}, & &\pinom{\K_i; c}{t}, & i\in \I, \, c &\in \Z, \, t\in \N_0
\end{align}
with slightly modified relations. Tagging them as in \cite{L-qgps-at-roots}, these are:
\begin{align}
\tag{b1}  &\text{the generators \eqref{eq:V0-gens-multlaced} commute with each other,} \label{eq:b1}
\\ \tag{b2} \K_i \K_i^{-1} & = 1, \quad \pinom{\K_i; c}{0} = 1, \label{eq:b2}
\\ \tag{b3} \label{eq:b3}
\pinom{\K_i; 0}{t}\pinom{\K_i; -t}{t'} &= {\pinom{t+t'}{t}}_{v^{d_i}} \pinom{\K_i; 0}{t + t'},\quad t,t' \geq 0,
\\ \tag{b4} \pinom{\K_i; c}{t} &- v^{-d_it} \pinom{\K_i; c+1}{t} = 
- v^{-d_i (c+1)} \K_i^{-1} \pinom{\K_i; c}{t - 1},
\quad t \geq 1, \label{eq:b4}
\\
\tag{b5} (v^{d_i} -v^{-d_i}) &\pinom{\K_i; 0}{1} = \K_i - \K_i^{-1}. \label{eq:b5}
\end{align}

The algebra $\V^0$ is related to $V^0$ in the following way.
For $i \in \I$, let $V_i^0$, respectively $\V_i^0$, be the subalgebra of $V^0$, respectively $\V^0$,
generated by $K_i^{\pm 1}$ and $\pinom{K_i; c}{t}$, respectively $\K_i^{\pm 1}$ and $\pinom{\K_i; c}{t}$,
$c \in \Z$, $t\in \N_0$. Then \eqref{eq:g6} and Proposition \ref{prop:V0-prop-lusztig}, respectively 
\eqref{eq:b1} and \cite[Theorem 6.7]{L-qgps-at-roots} imply that there are algebra isomorphisms
\begin{align}\label{eq:isom-Vi}
V^0 &\simeq V^0_1 \ot V^0_2 \ot \dots \ot V^0_{\theta}, &
\V^0 &\simeq \V^0_1 \ot \V^0_2 \ot \dots \ot \V^0_{\theta}.
\end{align}

\begin{lemma}\label{lemma:V0-isom-mult-laced}
	Let $\widetilde{\cA} = \cA$ regarded as $\cA$-algebra via $v \mapsto v^{d_i}$. Then
	$\V^0_i \simeq V^0_i \otimes_\cA \widetilde{\cA}$ as algebras.
\end{lemma}	

\pf First we claim that there is an algebra map $\psi_i: \V^0_i \to V^0_i \otimes_\cA \widetilde{\cA}$
given by $\K_i^{\pm 1} \mapsto K_i^{\pm 1} \ot 1$ and $\pinom{\K_i; c}{t} \mapsto \pinom{K_i; c}{t} \ot 1$.
Indeed, the images satisfy \eqref{eq:b2} by \eqref{eq:g7} and \eqref{eq:b5} by \eqref{eq:g5}.
Taking $c = t$ and $y = t'$ in \eqref{eq:g9} and inserting the right hand side in \eqref{eq:g9} we get 
\eqref{eq:b3}, while \eqref{eq:b4} follows applying \eqref{eq:g9} to both sides. The claim is proved 
and implies in turn that $\psi_i$ is an isomorphism, as it sends a basis to a basis
by Proposition \ref{prop:V0-prop-lusztig} and \cite[Theorem 6.7]{L-qgps-at-roots}.
\epf

From Lemma \ref{lem:comult-V0},
\eqref{eq:isom-Vi} and Lemma \ref{lemma:V0-isom-mult-laced} we see that $\V^0$ is a Hopf algebra
over $\cA$ with comultiplication determined by the $K_i$'s being group-likes. Let
\begin{gather*} 
\ku_{i,t} := \pinom{\K_i; 0}{t}, 
\\
\h_{i,n} := \frac{\K_i^{n}-\K_i^{-n}}{n(v^{d_in}-v^{-d_in})} \K_i^{n} = \frac{1}{n}\Big(\sum_{s \in \I_n} p_{n,s}(v^{d_i}) \, \ku_{i,s} \K_i^{\phi_n(s)} \Big) \K_i^{n} \in \V^0,
\end{gather*} 
$t, n\in \N$, $i\in\I$. Let $U(\hgo)\simeq \cB[t_i: i\in\I]$ as above and let 
$\Gamma = (\Z/{2\ell})^{\I}$ with generators  $(g_i)_{i\in \I}$. From the previous considerations we conclude:

\begin{prop}\label{th:V0-isom-mult-laced}
	Assume that $\xi^{2d_i} \neq 1$ for all $i\in \I$. 
	Then $V^0_{\cB} \simeq \cB\Gamma \otimes U(\hgo)$ as Hopf algebras via $g_i \mapsto \K_i$, $t_i  \mapsto \h_{i,\ell}$.
	\qed
\end{prop}


\begin{thebibliography}{AA}

\bibitem[A]{A-canad} N. Andruskiewitsch, \emph{Notes on extensions of Hopf algebras}. Canad. J. Math. \textbf{48} (1996), 3--42.


\bibitem[AA]{AA17} N. Andruskiewitsch, I. Angiono, 
\emph{On Finite dimensional Nichols algebras of diagonal type}. 
Bull. Math. Sci. \textbf{7} 353--573 (2017).
 
\bibitem[AAR1]{AAR1} N. Andruskiewitsch, I. Angiono, F. Rossi Bertone. \emph{The quantum divided power algebra of a finite-dimensional Nichols algebra of diagonal type}, Math. Res. Lett. \textbf{24} (2017), 619--643.

\bibitem[AAR2]{AAR2} \bysame \emph{A finite-dimensional Lie algebra arising from a Nichols algebra of   diagonal type (rank 2)}, Bull. Belg. Math. Soc. Simon Stevin 24 (1) (2017), 15--34.

\bibitem[AAR3]{AAR3} \bysame \emph{Lie algebras arising from Nichols algebras of diagonal type} \texttt{arXiv:1911.06586}.

\bibitem[AAY]{AAY} N. Andruskiewitsch, I. Angiono, M. Yakimov. \emph{Poisson orders on large quantum groups}, in preparation.


\bibitem[Ang1]{Ang-crelle} I. Angiono,
{\em On Nichols algebras of diagonal type}, J. Reine Angew. Math.  683 (2013), 189--251.

\bibitem[Ang2]{A-preNichols} \bysame \emph{Distinguished pre-Nichols algebras}, Transform. Groups \textbf{21} (2016), 1--33.

\bibitem[DK]{DK-root of 1} C. De Concini, V. Kac, \emph{Representations of quantum groups at roots of $1$}. Operator algebras, unitary representations, enveloping algebras, and invariant theory (Paris, 1989), 471--506, Progr. Math., \textbf{92} (1990), Birkhauser, Boston, MA.

\bibitem[DP]{DP-quantum} C. De Concini, C. Procesi, \emph{Quantum groups}. D-modules, representation theory, and quantum groups, 31--140, Lecture Notes in Math. \textbf{1565}, Springer, 1993.

\bibitem[H]{H-classif RS} I. Heckenberger, \emph{Classification of arithmetic root systems}, Adv. Math. \textbf{220} (2009) 59-124.



\bibitem[L]{Laugwitz} R. Laugwitz, 
\emph{Pointed Hopf algebras with triangular decomposition. A characterization of multiparameter quantum groups}, 
Algebr. Represent. Theory \textbf{19} 547--578 (2016).


\bibitem[Le1]{Lentner2} S. Lentner, \emph{A Frobenius homomorphism for Lusztig's quantum groups for arbitrary roots of unity}. Commun. Contemp. Math. 18 (2016), no. 3, 1550040, 42 pp.

\bibitem[Le2]{Lentner} \bysame, \emph{The unrolled quantum group inside Lusztig’s quantum group of divided powers}. 
Lett. Math. Phys.   \textbf{109} (2019), 1665--1682.



\bibitem[Lu1]{L-AiM} G. Lusztig, \emph{Quantum deformations of certain simple modules over enveloping algebras}, Adv. Math. 
\textbf{70} (1988), 237--249.

\bibitem[Lu2]{L-contemporary}  \bysame,   \emph{Modular representations and quantum groups}. Contemp. Math. \textbf{82} (1989),  59--77.  

\bibitem[Lu3]{L-fdHa-JAMS} \bysame, \emph{Finite dimensional Hopf algebras arising from quantized universal enveloping algebras}. J. Amer. Math. Soc. \textbf{3} (1990), 257--296.

\bibitem[Lu4]{L-qgps-at-roots}  \bysame,   \emph{Quantum groups  at  roots  of 1}. Geom. Dedicata \textbf{35} (1990),  89--114.  

\bibitem[Lu5]{L-libro} \bysame, \emph{ Introduction to quantum groups}. Birkh\"auser (1993).


\bibitem[Ma]{Majid} S. Majid, \emph{Double-bosonization of braided groups and the construction of $U_q(g)$}. Math. Proc. Camb.
Philos. Soc. \textbf{125}, 151--192 (1999).

\bibitem[M]{Mo-libro} S. Montgomery. 
\emph{Hopf algebras and their actions on rings}, CMBS \textbf{82},  
Amer. Math. Soc. (1993).

\bibitem[R]{Rad-libro} Radford, D. E., Hopf algebras, Series on Knots and Everything 49. 
Hackensack, NJ: World Scientific. xxii, 559 p.  (2012).


\bibitem[S]{S} Y. Sommerhauser, 
\emph{Deformed Enveloping Algebras}, 
New York J. Math. \textbf{2} (1996), 35--58.

\end{thebibliography}
\end{document}